\documentclass[9pt,journal]{IEEEtran}

\usepackage{amsfonts}
\usepackage{mathrsfs}
\usepackage{multicol}
\usepackage{multirow}
\usepackage{colortbl}
\usepackage{amssymb,amsmath,graphicx}

\usepackage{cite}

\newtheorem{Remark}{\bf Remark}[section]

\newenvironment{Proof}{\noindent{\em Proof:\/}}{\hfill $\Box$\par}
\newtheorem{Theorem}{\bf Theorem}[section]
\newtheorem{Lemma}{\bf Lemma}[section]
\newtheorem{Assumption}{\bf Assumption}[section]

\begin{document}

\title{Fixed-time Synchronization of Networked Uncertain Euler-Lagrange Systems}

\author{Yi~Dong and Zhiyong Chen
\thanks{This work has been supported in part by Shanghai Municipal Science and Technology Major Project under grant
2021SHZDZX0100, in part by National Natural Science
Foundation of China under grant 62073241
 and in part
by the Fundamental Research Funds for the Central Universities under grant 22120210127.}
\thanks{Y. Dong is with College of Electronic and Information Engineering,
Tongji University, Shanghai 200092, China. Email: yidong@tongji.edu.cn}
\thanks{Z. Chen is with  with the School of Electrical Engineering and Computing,
The University of Newcastle, Callaghan, NSW 2308, Australia. E-mail:
 zhiyong.chen@newcastle.edu.au}
 }

\maketitle

\begin{abstract}
This paper considers the  fixed-time control problem of a multi-agent system composed of a class of  Euler-Lagrange dynamics with parametric uncertainty and a dynamic leader under a directed communication network. A distributed   fixed-time  observer is first proposed to estimate the desired trajectory and then a  fixed-time controller is constructed by transforming uncertain Euler-Lagrange systems into second-order systems and utilizing the backstepping design procedure.
The overall design guarantees that the synchronization errors converge to zero in a prescribed time independent of initial conditions.
The control design conditions can also be relaxed for a weaker finite-time control requirement.
\end{abstract}


\IEEEpeerreviewmaketitle

\section{Introduction}


Fixed-time control for multi-agent systems, requiring exact achievement of a collective behavior
in a prescribed time independent of initial conditions, or finite-time control of a weaker requirement allowing
the prescribed time dependent on initial conditions, has attracted researchers' extensive attention over the past years due to its potential advantages in transient performance and robustness property \cite{bhat2000}.
The early work on finite-time formation control of single-integrator multi-agent systems can be found in \cite{xiao2009}.
For the leader-following consensus problem of general linear multi-agent systems, \cite{fu2016} proposed two classes of finite-time observers to estimate the second-order leader dynamics,  which can work in undirected and directed communication networks, respectively.
More efforts have also been devoted to  nonlinear systems. For example,
\cite{cao2014} considered the finite-time control of  first-order multi-agent systems with unknown nonlinear dynamics, while
both first-order and second-order nonlinear systems were considered in\cite{du2020}.  In particular,
observer-based control was proposed to solve the leader-following fixed-time consensus problem
under the strongly connected communication network.
 The fixed-time consensus problem was also investigated for double-integrator systems under directed communication network and
 more general multi-agent systems with high-order integrator dynamics in \cite{zuo2015,zuo2017}, respectively.

Euler-Lagrange systems capture a large class of contemporary engineering problems and
finite-time control of this class of systems has been intensively investigated, especially in
the individual setting.  For example,
\cite{hong2002} considered finite-time control for an Euler-Lagrange system based on the method for a double-integrator system, while \cite{hong2006,lin1} further dealt with nonlinear systems in the presence of uncertainties.
The work in \cite{yu2005} studied
a non-singular sliding surface and constructed a continuous finite-time control strategy for uncertain Euler-Lagrange system.
Furthermore, \cite{zhao2009}  designed an adaptive controller to track a desired trajectory in finite time and \cite{galicki2015}
proposed a method for handing both uncertain dynamics and globally unbounded disturbances.

The research on fixed-time or finite-time control of uncertain Euler-Lagrange systems in a network setting is relatively rare.
Some related results can be found in  \cite{huang2015} where, by  adaptive control technique, a finite-time synchronization controller was constructed
for a multi-agent system modeled by some mechanical nonlinear systems
 with a connected communication network.
 The recent work reported in \cite{hu2019finite} studied finite-time coordination behavior
 of a multiple Euler-Lagrange system with an undirected network in the absence of uncertainties.
 In particular, with the introduction of auxiliary variables,  the system can be converted into a simpler
 form such that the adding a power integrator method can be applied to ensure the convergence.

This paper provides a solution to the leader-following  fixed-time synchronization problem for multiple Euler-Lagrange systems with parametric uncertainty. The strategy is based on a class of observers that can accurately estimate a dynamic trajectory in a fixed time. The design  relaxes the undirected  and connected assumption for the communication network in \cite{du2020,zuo2017,huang2015,hu2019finite} and considers a directed network graph.  Then an observer-based controller is  proposed for the multi-agent system composed of a dynamic leader and multiple heterogeneous Euler-Lagrange dynamics,
as opposed to the finite-time control method for multiple special mechanical systems in \cite{huang2015}.
In particular, the distributed  control law is able to guarantee
each Euler-Lagrange system can track a desired trajectory in a prescribed time, independent of initial conditions.
 It is worth mentioning that the control design conditions can be relaxed for a weaker finite-time control requirement.
Also, a reduced continuous controller can be directly applied to the fixed-time synchronization problem
for second-order nonlinear systems with a directed graph.

Throughout the paper, we use the following notations.
For a vector $x =[x_1,\cdots,x_n]^T  \in\mathbb{R}^n$,
$\| x \|_1 =|x_1| +\cdots +|x_n|$  represents its Manhattan (${\cal L}_1$) norm,
$\| x \| =\sqrt{x_1^2+ \cdots, x_n^2}$ its Euclidean (${\cal L}_2$) norm,
and $|x| = [|x_1|,\cdots,|x_n|]^T$ its element-wise absolute valued vector.
For a matrix $X$, $|X|$ is also defined as its element-wise absolute valued matrix.
The power function operator is  element-wise in terms of $x^a =[x^a_1,\cdots,x^a_n]^T$ for $a>0$.
For two vectors (matrices) $X$ and $Y$, comparison operators are element-wise;
for example, $X\geq Y$ means $x_{ij} \geq y_{ij}$ for every
$x_{ij}$ and $y_{ij}$, the $(i,j)$-elements of $X$ and $Y$, respectively.
 The operator $\mbox{sig}^a(x)=[\mbox{sign}(x_1)|x_1|^a,\cdots,\mbox{sign}(x_n)|x_n|^a]^T$
 is defined for $a>0$ and the sign function $\mbox{sign}(\cdot)$.

\section{Problem Formulation}
Consider a group of  $m$-link robotic manipulators of the following Euler-Lagrange dynamics
  \begin{align}
    M_i(q_i)\ddot{q}_i+ C_i(q_i,\dot{q}_i)\dot{q}_i+G_i(q_i)=\tau_i,\;i=1,\cdots,N, \label{1}
  \end{align}
where $q_i\in \mathbb{R}^m$, $\dot{q}_i\in \mathbb{R}^m$ are the vectors of generalized
position and velocity of the $i$-th  robotic manipulator, also called agent $i$,
$M_i(q_i)\in \mathbb{R}^{m \times m}$ is a symmetric and
positive definite inertia matrix, $C_i(q_i,\dot{q}_i)\dot{q}_i\in
\mathbb{R}^{m}$ contains the Coriolis and centrifugal forces, $G_i(q_i)\in \mathbb{R}^m$ is the gravitational torque, and $\tau_i\in \mathbb{R}^m$ is the vector of control force.
 The reference is generated by a leader system, called agent 0,  described  as follows:
\begin{align}
  \dot{\eta}_0&=S\eta_0,\; q_0=E\eta_0, \label{2}
\end{align}
where $\eta_0\in\mathbb{R}^n$ is the state, $q_0\in \mathbb{R}^m$ is the desired trajectory to track,
and  $S\in \mathbb{R}^{n\times n}$,  $E\in\mathbb{R}^{m\times n}$ are constant matrices.

The multi-agent system under consideration is composed of the $N$ dynamics in (\ref{1}) and the dynamic leader (\ref{2}).
The information flow among all the $N+1$ agents is described by a digraph $\mathcal{G}=(\mathcal{V}, \mathcal{E})$ where
$\mathcal{V}=\{0,1,\cdots,N\}$ is  the node set and $\mathcal{E} \subset \mathcal{V}\times \mathcal{V}$ is the edge set.
Each element $(j,i)\in\mathcal{E}$ represents
the edge from agent $j$ to agent $i$.
For $i,j \in \mathcal{V}$, $a_{ii}=0$, $a_{ij}>0$ if $(j,i)\in\mathcal{E}$ and $a_{ij}=0$ otherwise.
Let $H=[h_{ij}]_{i,j=1}^N$ be the Laplacian matrix of the subnetwork composed of agents $1,\cdots, N$,
where $h_{ii}=\sum_{j=0}^N a_{ij}$ and $h_{ij}=-a_{ij}$  for $i\neq j$, $i,j=1,\cdots,N$.


The objective of this paper is to design a distributed control law such that each agent of (\ref{1}) can track the desired trajectory $q_0$ in   fixed time. More specifically, we
consider the class of control laws of the form
  \begin{align}
   \tau_i&=f_{1i}(q_i,\dot{q}_i,\eta_i,\dot{\eta}_i) \nonumber\\
   \dot{\eta}_i&=f_{2i}(\eta_i,\sum_{j=0}^Na_{ij} (\eta_i-\eta_j)), \; i=1,\cdots,N. \label{mlaw2}
  \end{align}
Let
\begin{align*}
  x_i=\left[ \begin{array}{c}
q_i-q_0\\
\dot{q}_i-\dot{q}_0\\
\eta_i-\eta_0
\end{array}
\right],i=1,\cdots,N,\;
x=\left[ \begin{array}{c}
x_1\\
\vdots\\
x_N
\end{array}
\right]\in\mathbb{R}^{n_x}.
\end{align*}
Suppose the closed-loop system composed of (\ref{1}), (\ref{2}) and (\ref{mlaw2})
possesses unique solutions in forward time for all initial conditions.
 Then the fixed-time synchronization problem can be described as follows  based on
the concept of fixed-time stability \cite{polyakov2011}.

 {\bf Fixed-time synchronization  problem:}  Given the system composed of
 (\ref{1}) and  (\ref{2}) with  the corresponding
digraph  $\mathcal{G}$, design a distributed control law of the form (\ref{mlaw2})
such that, for all initial conditions $x(0) =x_0$ and $\eta_0(0) =\eta_{00}$,
the equilibrium $x=0$ of the closed-loop system is  (globally) fixed-time stable.
That is, the solution $x(t)$  exists for $t\geq 0$ and $x=0$ is Lyapunov stable, and moreover,
there exists a fixed time $T^{*}$, independent of $x_0$ or $\eta_{00}$, such that
\begin{align}
& \lim_{t\rightarrow T^{*}} x(t)=0, \nonumber\\
& x(t)=0,\; t\geq T^{*},\; \forall  x_0 \in \mathbb{R}^{n_x},  \eta_{00} \in \mathbb{R}^{n}. \label{obj}
\end{align}

\begin{Remark}
 When the existence of a fixed time $T^{*}$ is relaxed to the existence
of a settling-time function $T(x_0,\eta_{00})$, the above fixed-time synchronization problem
is called a finite-time synchronization problem, which is based on the definition of finite-time stability \cite{bhat2000}.
In practice, there is similarity between finite-time convergence and asymptotical (exponential)
convergence, both of which require the  convergence of trajectories to (proximity of) an equilibrium point in a
finite amount of time which depends on the initial conditions. But
fixed-time convergence is a more practically interesting feature which requests
it happen in a prescribed time independent of initial conditions.
There are more constraints on the controller design conditions that will be studied in this paper.
\end{Remark}

For the solvability of the aforementioned problem, we need the following standard assumption on  the communication network.
\begin{Assumption}\label{a1}
The graph $\mathcal{G}$  contains a spanning tree with node 0 as the root.
\end{Assumption}

\begin{Remark}
Under Assumption \ref{a1}, all the eigenvalues of $H$ have positive real parts; see, e.g.,  \cite{su2015}. By Theorem 2.5.3 of \cite{horn1994},  there exists a positive definite diagonal  matrix $\bar D\in\mathbb{R}^{N\times N}$ such that $H^T \bar D+ \bar DH$ is  positive definite. Let $\lambda_m >0$ be the smallest eigenvalue of $H^T \bar D+ \bar DH$ and
$D=\mbox{diag}(d_1,\cdots,d_N) = 2 \bar D / \lambda_m$. One has $H^TD+DH \geq2I_N$.
\end{Remark}


We end this section with some technical lemmas from, e.g., \cite{zuo2015}, \cite{polyakov2011}, \cite{hardy1952} and \cite{qian2001},
which will be used in the proofs of the main results in this paper.

\begin{Lemma} \label{ha}
For any $\xi_i\in\mathbb{R}$, $i=1,\cdots,n$,  and any $p\in(0,1]$,  $(\sum_{i=1}^n|\xi_i|)^p\leq \sum_{i=1}^n |\xi_i|^p\leq n^{1-p} (\sum_{i=1}^n|\xi_i|)^p$ \cite{hardy1952}. For any $p>1$, $\sum_{i=1}^n|\xi_i|^p\leq (\sum_{i=1}^n |\xi_i|)^p\leq n^{p-1} \sum_{i=1}^n|\xi_i|^p$ \cite{zuo2015}.
\end{Lemma}

\begin{Lemma}  \cite{hardy1952} \label{ha1}
The inequality $|\xi_i^p-\xi_j^p|\leq 2^{1-p}|\xi_i-\xi_j|^p$ holds for $\forall  \xi_i,\xi_j\in\mathbb{R}$ and $0<p\leq 1$ and $p$ is a ratio of two odd integers.
\end{Lemma}

\begin{Lemma}  \cite{qian2001} \label{qi}
 The inequality $|\xi_i|^c|\xi_j|^d\leq \frac{c}{c+d}r|\xi_i|^{c+d}+\frac{d}{c+d}r^{-\frac{c}{d}}|\xi_j|^{c+d}$ holds for $\forall \xi_i,\xi_j\in\mathbb{R}$ and $c,d,r>0$.
\end{Lemma}


\begin{Lemma} (Lemma 1, \cite{polyakov2011}) \label{poly}
Consider the system $\dot z = \phi (z,t)$ where $\phi: \mathbb{R}^l \times [0,\infty) \mapsto \mathbb{R}^l$
satisfies $\phi(0,t)=0$.
Suppose there exits a continuously differentiable function $V:~ \mathbb{R}^l \mapsto \mathbb{R}$ such that
(i) $V$ is positive  definite and proper; and (ii) there exist real numbers
$p_0,q_0,p,q,k>0$ with  $pk<1$ and $qk>1$ such that
$\dot V(z) \leq - (p_0 (V(z))^{p}+q_0 (V(z))^q )^k$.
 Then, the equilibrium $z=0$ is (globally) fixed-time stable and there is a constant settling-time
   $T^{*}\leq \frac{1}{p_0^k(1-pk)}+\frac{1}{q_0^k(qk-1)}$.
\end{Lemma}

\section{Distributed observer design}

As the agents not  connected to agent 0 do not have access to the information of the dynamic leader (\ref{2}),
its state needs to be estimated by a properly designed fixed-time observer as follows:
\begin{align}\label{ob}
  \dot{\eta}_i&=S\eta_i-c_1 y_i-c_2 \mbox{sig}^{a}(y_i)-c_3 \mbox{sig}^{b}(y_i),   \nonumber\\
  y_i & =\sum_{j=0}^Na_{ij} (\eta_i-\eta_j),\;   i=1,\cdots,N.
\end{align}
 In this section, we construct a lemma  based on
the fixed-time observer (\ref{ob}).
Let $\eta^T = [\eta^T_0, \eta^T_1, \cdots, \eta^T_N ]^T$
for the convenience of presentation.

\begin{Lemma}  \label{lemma_fixobs}
Consider the system composed of (\ref{2}) and (\ref{ob}) under Assumption \ref{a1} with
$0<a<1$, $b>\frac{1}{a}>1$, $c_1>\|D\otimes S\|$ and $c_2,c_3>0$.
   There exists a constant settling-time $T^*_1\geq 0$ such that,
$\forall \eta(0) \in \mathbb{R}^{(N+1)n}$,
\begin{align}
& \lim_{t\rightarrow T_1^*  }(\eta_i(t)-\eta_0(t))=0, \nonumber\\
& \eta_i(t)-\eta_0(t)=0,\; t\geq T_1^*, \; i= 1,\cdots, N. \label{T1star}
\end{align}
\end{Lemma}

\begin{Proof}
Let $\bar{\eta}_i=\eta_i-\eta_0$, $i=0,1,\cdots,N$. The observer (\ref{ob}) can be rewritten as
\begin{align}\label{obb}
  \dot{\bar{\eta}}_i &=S\bar{\eta}_i-c_1 y_i-c_2 \mbox{sig}^{a}(y_i)-c_3 \mbox{sig}^{b}(y_i),\nonumber\\
  y_i&=\sum_{j=0}^N a_{ij} (\bar{\eta}_i-\bar{\eta}_j),\; i=1,\cdots,N.
\end{align}
Let $Y_i=c_1 y_i+c_2\mbox{sig}^{a}(y_i)+c_3 \mbox{sig}^{b}(y_i)$ and $\bar{\eta}$, $y$, $Y$ be the column stacks of $\bar{\eta}_i$, $y_i$, $Y_i$, $i=1,\cdots,N$.
Note $y=(H\otimes I_n) \bar{\eta}$ and
\begin{equation}\label{yc}
  \dot{y}=(I_N\otimes S)y-(H\otimes I_n)Y.
\end{equation}
And, let  \begin{align}
  V(y)=& \sum_{i=1}^N \left( \frac{c_2 d_i }{1+a}\|y_{i}^{1+a}\|_1+ \frac{c_3 d_i}{1+b}\|y_{i}^{1+b}\|_1 \right) \nonumber\\
  &+\frac{c_1}{2}y^T(D\otimes I_n)y. \label{VT}
\end{align}
Along the trajectory of (\ref{yc}), the time derivative of $V(y)$ satisfies
  \begin{align*}
     \dot{V}(y)
    = &\sum_{i=1}^N d_i (c_2   \mbox{sig}^a (y_{i})+c_3   \mbox{sig}^b (y_{i}) )^T  \dot{y}_{i}
    +c_1y^T(D\otimes I_n)\dot{y}\\
  = & Y^T(D\otimes I_n)\dot{y}\\
   = &Y^T(D\otimes  S)y-\frac{1}{2}Y^T((H^TD+DH)\otimes I_n)Y\\
  \leq  & \frac{1}{2}Y^TY+\frac{1}{2} \|D\otimes  S \|^2y^Ty-Y^TY\\
    = & -\frac{1}{2} (\|Y\|^2 -\|D\otimes  S\|^2 \|y\|^2).
  \end{align*}
Further calculation shows that \begin{align*}
    \|Y_i\|^2 =&\|c_1 y_i+c_2\mbox{sig}^{a}(y_i)+c_3 \mbox{sig}^{b}(y_i))\|^2  \\
    \geq &c_1^2 \|y_i\|^2 +c_2^2 \| \mbox{sig}^{a}(y_i)\|^2  +c_3^2 \| \mbox{sig}^{b}(y_i)\|^2 \\
   \geq &c_1^2 \|y_i\|^2 +c_2^2 \|y_i^{2a}\|_1  +c_3^2 \|y_i^{2b}\|_1 .
        \end{align*}
 By Lemma~\ref{ha}, for $0<a<1$ and $b>1$,
\begin{align*}
&\sum_{i=1}^N \|y_i^{2a}\|_1 \geq(\sum_{i=1}^N \|y_i^{2}\|_1 )^a =(\|y\|^2 )^a,\\
&\sum_{i=1}^N \|y_i^{2b}\|_1 \geq  \frac{1}{(nN)^{b-1}} (\|y\|^2 )^b.
        \end{align*}
 As a result,
  \begin{align*}
    \|Y\|^2 =  &  \sum_{i=1}^N  \|Y_i\|^2
    \geq
    c_1^2 \|y\|^2+c_2^2(\|y\|^2)^a+\frac{c_3^2}{(nN)^{b-1}}(\|y\|^2)^b
  \end{align*}
and hence
  \begin{align}
  \dot{V}(y) \leq & -\frac{1}{2}(c_1^2-\|D\otimes S\|^2) \|y\|^2 -\frac{c_2^2}{2} (\|y\|^2)^a \nonumber\\
&- \frac{c_3^2}{2(nN)^{b-1}}(\|y\|^2)^b \nonumber\\
\leq & -\hat{c}_1\left( \|y\|^2+(\|y\|^2)^a+(\|y\|^2)^b \right) \label{DVR}
  \end{align}
for
   \begin{align*}
\hat{c}_1=\min\left\{\frac{1}{2}(c_1^2-\|D\otimes  S\|^2),\frac{c_2^2}{2},\frac{c_3^2}{2(nN)^{b-1}}\right\} >0.
\end{align*}

Analysis on  \eqref{VT} using Lemma~\ref{ha} and noting
$0<\frac{2a}{1+a}<1$ and $\frac{a(1+b)}{1+a}>1$ gives
\begin{align*}
V^{\frac{2a}{1+a}} \leq& \sum_{i=1}^N
\Big((\frac{c_1d_i}{2})^{\frac{2a}{1+a}} \|y_{i}^{\frac{4a}{1+a}}\|_1 +(\frac{c_2 d_i}{1+a})^{\frac{2a}{1+a}} \|y_{i}^{2a} \|_1\\
 &+(\frac{c_3 d_i}{1+b})^{\frac{2a}{1+a}}\|y_{i}^{\frac{2a(1+b)}{1+a}}\|_1 \Big)\\
 \leq& \hat{c}_2
 \left( (\|y\|^2)^{\frac{2a}{1+a}}+ (\|y\|^2)^a+(\|y\|^2)^{\frac{a(1+b)}{1+a}} \right)
\end{align*}
for $d_{\max}= \max\{d_1,\cdots,d_N\}$ and
 \begin{align*}
\hat{c}_2=& \max\Big\{(\frac{c_1d_{\max}}{2})^{\frac{2a}{1+a}} (nN)^{1-\frac{2a}{1+a}},\\
&(\frac{c_2 d_{\max}}{1+a})^{\frac{2a}{1+a}}(nN)^{1-a}, (\frac{c_3d_{\max}}{1+b})^{\frac{2a}{1+a}} \Big\}.
\end{align*}
Since $a<\frac{2a}{1+a}<1$ and $a<\frac{a(1+b)}{1+a}<b$, we can easily verify
\begin{align*}
(\|y\|^2)^{\frac{2a}{1+a}}\leq &\|y\|^2 +(\|y\|^2)^{a} \\
(\|y\|^2)^{\frac{a(1+b)}{1+a}} \leq& (\|y\|^2)^{a}  +(\|y\|^2)^{b}
\end{align*}
and hence
\begin{align}
V^{\frac{2a}{1+a}}   \leq& \hat{c}_2
 \left(\|y\|^2 +3 (\|y\|^2)^a+(\|y\|^2)^b \right). \label{V1}
\end{align}

Similarly, for $b>1$, $\frac{2b}{1+b}>1$ and $\frac{b(1+a)}{1+b}>1$,
\begin{equation*}
\begin{aligned}
 V^{\frac{2b}{1+b}}
 &\leq \hat{c}_3 \left((\|y\|^2)^{\frac{2b}{1+b}}+ (\|y\|^2)^{\frac{b(1+a)}{1+b}}+ (\|y\|^2)^b \right)
\end{aligned}
\end{equation*}
for
 \begin{align*}
\hat{c}_3= &\max\left\{(\frac{c_1d_{\max}}{2})^{\frac{2b}{1+b}},~
 (\frac{c_2 d_{\max}}{1+a} )^{\frac{2b}{b+1}},
 (\frac{c_3 d_{\max}}{1+b} )^{\frac{2b}{b+1}} \right\} \\
& \times (3nN)^{\frac{b-1}{b+1}}.
 \end{align*}
 Since $1<\frac{2b}{1+b}<b$ and $a<\frac{b(1+a)}{1+b}<b$, we can easily verify
\begin{align*}
(\|y\|^2)^{\frac{2b}{1+b}} \leq &\|y\|^2 +(\|y\|^2)^{b} \\
(\|y\|^2)^{\frac{b(1+a)}{1+b}} \leq& (\|y\|^2)^{a}  +(\|y\|^2)^{b}
\end{align*}
and hence
\begin{align}\label{V2}
 V^{\frac{2b}{1+b}}   \leq& \hat{c}_3
 \left(\|y\|^2 + (\|y\|^2)^a+3(\|y\|^2)^b \right)
\end{align}

 Finally,  by \eqref{V1} and \eqref{V2}, one has
\begin{align*}
\frac{1}{ \hat{c}_2} V^{\frac{2a}{1+a}}
+\frac{1}{ \hat{c}_3} V^{\frac{2b}{1+b}}
  \leq  4 \left( \|y\|^2 + (\|y\|^2)^a+ (\|y\|^2)^b \right), \end{align*}
which, compared with  \eqref{DVR}, implies
  \begin{align*}
  \dot{V}  \leq -\frac{\hat{c}_1}{4\hat c_2} V^{\frac{2a}{1+a}}
  -\frac{\hat{c}_1}{4\hat c_3}  V^{\frac{2b}{1+b}}.
  \end{align*}
 By Lemma~\ref{poly}, the system (\ref{yc}) is fixed-time stable.
 In particular, there exists a constant
 \begin{align*}
 T_1^{*}\leq \frac{4\hat{c}_2 (a+1)}{\hat{c}_1(1-a)}+\frac{4\hat{c}_3(b+1)}{\hat{c}_1(b-1)},
 \end{align*}
  such that
   $\lim_{t\rightarrow T_1^{*}} y(t)=0$ and
$y(t)=0$, $ t\geq  T_1^{*}$.
Under Assumption \ref{a1},  we have $\bar{\eta}=(H^{-1}\otimes I_n)y$ and hence  \eqref{T1star}.
The proof is thus completed.
\end{Proof}

\begin{Remark} When $c_3=0$, the observer \eqref{ob} reduces to a finite-time observer
\begin{align}\label{ob1}
  \dot{\eta}_i &=S\eta_i-c_1 y_i-c_2 \mbox{sig}^{a}(y_i), \nonumber\\
  y_i & =\sum_{j=0}^Na_{ij} (\eta_i-\eta_j),\;   i=1,\cdots,N.
\end{align}
Consider the system composed of (\ref{2}) and (\ref{ob1}) under Assumption~\ref{a1} with
$0<a<1$, $c_1>\| D\otimes S\|$ and $c_2>0$.
There exists a settling-time function $T_1(\eta(0)) \geq 0$ such that,  $\forall \eta(0) \in \mathbb{R}^{(N+1)n}$,
\begin{align}
& \lim_{t\rightarrow T_1  }(\eta_i(t)-\eta_0(t))=0, \; i= 1,\cdots, N\nonumber\\
& \eta_i(t)-\eta_0(t)=0,\; t\geq T_1(\eta(0)).  \label{T1}
\end{align}
The proof of the above statement follows that of Lemma~\ref{lemma_fixobs}
using simple arguments.  In particular, we can obtain
   \begin{align*}
  \dot{V}  \leq  -\frac{\hat{c}_1}{3\hat c_2} V^{\frac{2a}{1+a}}.    \end{align*}
 In other words,  $\dot{V}  + \varrho  V^{\frac{2a}{1+a}}$ is negative definite
 for any $\varrho < \frac{\hat{c}_1}{3\hat c_2}$.
By Theorem 1 in \cite{bhat2000},  the system (\ref{yc}) is finite-time stable.
 In particular, there exists a finite settling-time function
 \begin{align*} \bar T_1(y(0))\leq \frac{3\hat{c}_2 (a+1) V(y(0))^{1-\frac{2a}{1+a}}}{\hat{c}_1(1-a)},
 \end{align*}
  such that
   $\lim_{t\rightarrow \bar T_1(y(0)) } y(t)=0$ and
$y(t)=0$, $ t> \bar T_1(y(0))$, $\forall y(0) \in \mathbb{R}^{Nn}$.
Under Assumption \ref{a1},  we have $\bar{\eta}=(H^{-1}\otimes I_n)y$ and hence  \eqref{T1}
for $T_1(\eta(0) ) = \bar T_1((H\otimes I_n)\bar{\eta}(0))$.
\end{Remark}

\section{ Robust Controller design}
Based on  the fixed-time observer (\ref{ob}),
we further propose a distributed robust control law to solve the leader-following fixed-time synchronization problem for
the multiple Euler-Lagrange systems.
 It is assumed  the model \eqref{1} contains uncertainties and the
terms $M_i(q_i)$, $C_i(q_i,\dot{q}_i)$, and $G_i(q_i)$ are not completely known, but they
satisfy the following bounded conditions
 \begin{align}
& k_{\underline{m}}I_m\leq M_i(q_i)\leq k_{\overline{m}}I_m , \nonumber\\
& \|C_i(q_i,\dot{q}_i) \|\leq k_c \|\dot{q}_i \| , \; \|G_i(q_i) \|\leq k_g ,\; \forall  q_i, \dot{q}_i \in \mathbb{R}^m, \label{km}
  \end{align}
for some positive constants $k_{\underline{m}}$, $k_{\overline{m}}$, $k_c$ and $k_g$.
Throughout the section, we consider every individual agent $i=1,\cdots,N$.

First, the equations in (\ref{1}) can be rewritten as, with $v_i = \dot{q}_i$,
  \begin{align*}
    \dot{q}_i=v_i,~\dot{v}_i=M_i^{-1}(q_i)(\tau_i-C_i(q_i,v_i)v_i-G_i(q_i)),
      \end{align*}
which is a second-order system in the presence
parametric uncertainty, i.e., the terms $M_i(q_i)$, $C_i(q_i,v_i)$ and $G_i(q_i)$ are unknown for  controller design.
Therefore,  the conventional  fixed-time control laws for second-order systems cannot be directly applied.
To introduce a new method,  we perform the following transformation:
  \begin{align*}
 \bar{q}_i=q_i-E\eta_i, \bar{v}_i=v_i-E\dot{\eta}_i, \tau_i=\hat{M} u_i, \;\hat{M}=\frac{2I_m}{k_{\underline{m}}^{-1}+k_{\overline{m}}^{-1}}.
   \end{align*}
Also, let $u_i=u_{1i}+u_{2i}$ with $u_{1i}$ and $u_{2i}$ to be designed.
As a result, the above equations become
  \begin{align*}
    \dot{\bar{q}}_i =  \bar{v}_i,~ \dot{\bar{v}}_i =&M_i^{-1}(q_i)\tau_i +F_i(q_i,v_i)-E\ddot{\eta}_i\\
    =&u_{1i}+u_{2i}+(M_i^{-1}(q_i)\hat{M}-I_m)(u_{1i}+u_{2i})\\
    &+F_i(q_i,v_i)-E\ddot{\eta}_i
  \end{align*}
where $F_i(q_i,v_i)=-M_i^{-1}(q_i)(C_i(q_i,v_i)v_i+G_i(q_i))$.
Moreover, it can be put in a compact form
\begin{equation}\label{dou}
  \dot{\bar{q}}_i=\bar{v}_i,~~\dot{\bar{v}}_i=u_{2i}-E\ddot{\eta}_i+Z_i
\end{equation}
with
  \begin{align}
 Z_i=u_{1i}+(M_i^{-1}(q_i)\hat{M}-I_m)(u_{1i}+u_{2i})+F_i(q_i,v_i). \label{Zi}
    \end{align}

Inspired by \cite{galicki2015}, we construct a lemma that
motivates the design of $u_{1i}$.

\begin{Lemma}\label{lemma_u1i} Consider the quantity $Z_i$ defined in \eqref{Zi} with
the control law
  \begin{align}
  &u_{1i}=\left\{
           \begin{array}{ll}
              -\frac{\kappa}{1-\epsilon}\frac{\zeta_i}{\|\zeta_i\|}(\epsilon \|u_{2i}\|+f_i(v_i)),& \|\zeta_i\|\neq 0\\
             0,&\|\zeta_i\|=0\\
           \end{array}
         \right.  \nonumber \\
&\kappa\geq 1, \; \epsilon=\frac{k_{\underline{m}}^{-1}-k_{\overline{m}}^{-1}}{k_{\underline{m}}^{-1}+k_{\overline{m}}^{-1}},\; f_i(v_i)=k_{\underline{m}}^{-1}(k_c \|v_i \|^2+k_g), \label{u1i}
  \end{align}
for an arbitrary $\zeta_i\in\mathbb{R}^m$. Then,
$\zeta_i^TZ_i \leq 0$ holds for any $u_{2i}\in\mathbb{R}^m$.
\end{Lemma}
\begin{Proof}
From the properties of Euler-Lagrange system, we have the following facts:
  \begin{align*}
 \|M_i^{-1}(q_i)\hat{M}-I_m\| = & \|\frac{2M_i^{-1}}{k_{\underline{m}}^{-1}+k_{\overline{m}}^{-1}}-I_m\|
    \leq
  \epsilon \\
 \|F_i(q_i,v_i)\|\leq& k_{\underline{m}}^{-1}(k_c\|v_i\|^2+k_g)= f_i(v_i),
   \end{align*}
   which will be used in the calculation below.

 When $\|\zeta_i\|=0$,  $\zeta_i^TZ_i(t)\leq 0$ holds trivially. Otherwise, one has the following direct calculation
  \begin{align*}
    \zeta_i^T Z_i \leq & \zeta_i^Tu_{1i} +\|\zeta_i \|(\|M_i^{-1}(q_i)\hat{M} -I_m\| (\|u_{1i}\|+\|u_{2i}\|)\\
&+\|F_i(q_i,v_i)\|)\\
\leq & \zeta_i^Tu_{1i} +\|\zeta_i\|(\epsilon (\|u_{1i}\|+\|u_{2i}\|)+f_i(v_i))\\
\leq & (-\frac{\kappa}{1-\epsilon}+1+\frac{\epsilon\kappa}{1-\epsilon})\|\zeta_i\|(\epsilon \|u_{2i}\|+f_i(v_i))\\
= &-(\kappa-1)\|\zeta_i\|(\epsilon \|u_{2i}\|+f_i(v_i))\leq 0,
  \end{align*}
which completes the proof.
  \end{Proof}

\begin{Remark}  As the system dynamics considered in this paper contain uncertainties, a robust control
approach is used in the design of $u_{1i}$. In particular, to guarantee $\zeta_i^TZ_i(t)\leq 0$, which will be used
later for proof of convergence, $u_{1i}$ is designed based on the boundaries of the uncertainties characterized
by \eqref{km} via high gain domination. It is worth noting that the control gains in $u_{1i}$ become higher
if  $k_c$ and $k_g$ are larger and/or  $\epsilon$ is closer to 1 (corresponding to a bigger difference between $k_{\underline{m}}$ and $k_{\overline{m}}$), i.e., the size of uncertainties is larger.
It is  a general principle in robust control that control gains depend on the size of uncertainties.
In practice, when system parameters cannot be precisely measured, a smaller range of uncertainties would be beneficial
for controller design.
\end{Remark}

With Lemma~\ref{lemma_u1i}  ready for  $u_{1i}$, the remaining task is to select a specific  $\zeta_i$
and design  $u_{2i}$  such that the second-order system  (\ref{dou}) is fixed-time stable, which is more complicated than finite-time control; see some existing methods in \cite{du2020,zuo2015},\cite{song2017}.
For solving such problem, we first introduce an explicit procedure of designing a set of parameters to be used for the controller design.
It is worth noting  that these parameters are independent of system dynamics.
Let $\frac{1}{2}<\alpha<1$ and $\beta>1$ be two specified rational numbers of ratio of two odd integers.
Define four constants
  \begin{align*}
  p_1  =0, \; p_2 =\beta-\alpha, \;
  p_3 ={\frac{\beta}{\alpha}-\beta+\alpha-1},\; p_4={\frac{\beta}{\alpha}-1}
   \end{align*}
 and four functions, for  $p\geq 0$ and $\lambda>0$,
 \begin{align*}
l_1(p) &= \frac{2^{1-\alpha}p}{p+1+\alpha}+\frac{p+\alpha}{p+\alpha+1},~l_2(p)=\frac{p+\beta}{p+\beta+1}, \\
 l_3(p,\lambda)&=\frac{2^{1-\alpha}(1+\alpha)}{p+1+\alpha}\lambda^{\frac{p+\alpha+1}{1+\alpha}}+\frac{(\lambda\gamma_1)^{p+\alpha+1}}{p+\alpha+1},\\
l_4(p,\lambda)&=\frac{(\lambda\gamma_2)^{p+\beta+1}}{p+\beta+1}.
\end{align*}
For the convenience of presentation, we also define
\begin{align*}
\ell_1(p) &= (2-\alpha)2^{1-\alpha} l_1(p),~\ell_2(p) = (2-\alpha)2^{1-\alpha} l_2(p),\\
\ell_3(p,\lambda) &= (2-\alpha)2^{1-\alpha} l_3(p,\lambda),~\ell_4(p,\lambda) = (2-\alpha)2^{1-\alpha} l_4(p,\lambda).
\end{align*}
Next, pick
$L_1= \max\{ \ell_2(p_2),  \ell_1(p_3)  \}$
and two positive parameters
  \begin{align}
\gamma_1 > & \max \Big\{ \frac{2^{1-\alpha}}{1+\alpha} + \ell_1(p_1)     +2L_1  , \frac{2^{1-\alpha}\frac{\beta}{\alpha}}{\frac{\beta}{\alpha}+\alpha} + \ell_2(p_3) +\ell_1(p_4) \Big\} \nonumber\\
  \gamma_2  > & \max \left\{ \ell_2(p_4) +2L_1, \ell_2(p_1) +\ell_1(p_2)\right\} .\label{gamma1}
\end{align}
Now, it is ready to select
\begin{align*}
  \lambda_1=  \gamma_1^{\frac{1}{\alpha}},\;  \lambda_2= \gamma_1^{\frac{1}{\alpha}-1}\frac{\gamma_2\beta}{\alpha},  \;
  \lambda_3 = \gamma_1\gamma_2^{\frac{1}{\alpha}-1} ,\;
  \lambda_4 = \gamma_2^{\frac{1}{\alpha}}  \frac{\beta}{\alpha}.
   \end{align*}
Then, pick
 \begin{align*}
L_2 = &\max\Big\{ \frac{ 2^{1-\alpha}\alpha}{\frac{\beta}{\alpha}+\alpha}+ \ell_4(p_3,\lambda_3)+ \ell_3(p_4, \lambda_4),\\
& \ell_4(p_1,\lambda_1) +\ell_3(p_2, \lambda_2) ,  \ell_4(p_2, \lambda_2),\ell_3(p_3, \lambda_3)    \Big\}.
\end{align*}
Finally, we select the following two parameters
  \begin{align}
k_1   > \frac{2^{1-\alpha}\alpha}{\alpha+1} + \ell_3(p_1, \lambda_1)    +4L_2,\;
k_2   > \ell_4(p_4, \lambda_4) + 4L_2.\label{kk}
\end{align}

With these parameters obtained, it is ready to have the following lemma.
\begin{Lemma} \label{unii}
Consider the system (\ref{dou}) where $u_{1i}$ is given in \eqref{u1i} with
  \begin{align} \label{zetai}
\zeta_i= \varepsilon_{i}^{2-\alpha}  \end{align}
and
 \begin{align}
  u_{2i}&=-k_1 \varepsilon_{i}^{2\alpha-1}-k_2 \varepsilon_{i}^{\frac{\beta}{\alpha}+\beta+\alpha-2}+ESS\eta_i, \nonumber\\
  \varepsilon_{i}&=\bar{v}_{i}^{\frac{1}{\alpha}} + (\gamma_1\bar{q}_{i}^{\alpha} + \gamma_2\bar{q}_{i}^{\beta})^{\frac{1}{\alpha}}.
  \label{unfix}
\end{align}
Suppose the observer governing $\eta_i$ satisfies Lemma~\ref{lemma_fixobs}.
If the control parameters $\gamma_1,\gamma_2,k_1,k_2$
satisfy (\ref{gamma1}) and (\ref{kk}), then the equilibrium of  (\ref{dou}) at the origin is  fixed-time stable.
In particular, there exists a constant settling-time  $T_{2}^* \geq  0$ such that
\begin{align}
& \lim_{t\rightarrow T_1^*+T^*_{2}  } [\bar q_i(t), \bar v_i(t)] =0, \nonumber\\
&  [\bar q_i(t), \bar v_i(t)]=0,\; t\geq T_1^* + T^*_{2} , \;
  \forall  \bar q_i(T^*_1), \bar v_i(T^*_1) \in \mathbb{R}^{m}.  \label{T2star}
\end{align}

\end{Lemma}

\begin{Proof} For the convenience of proof,
we define the following variables
\begin{align*}
 \bar{\varepsilon}_i &= \varepsilon_{i}^{2\alpha-1},~~ \hat{\varepsilon}_i= \varepsilon_{i}^{\frac{\beta}{\alpha}+\beta+\alpha-2}\\
 \hat{v}_{i}^{*}&=-\gamma_1\bar{q}_{i}^{\alpha}-\gamma_2\bar{q}_{i}^{\beta},~~\varepsilon_{i} =\bar{v}_{i}^{\frac{1}{\alpha}}-\hat{v}_{i}^{*\frac{1}{\alpha}}
\end{align*}
Let
 \begin{align*}
W_{1i} (\bar q_i) =\frac{1}{2}\|\bar{q}_{i}\|^2+\frac{1}{\beta/\alpha+1} \| \bar{q}_{i}^{\frac{\beta}{\alpha}+1}\|_1.
 \end{align*}
 Along the trajectory of $\bar{q}_i-$th subsystem in (\ref{dou}), the derivative of $W_{1i} (\bar q_i)$ satisfies
  \begin{align*}
    \dot W_{1i} (\bar q_i)
    = & (\bar{q}_i +\bar{q}_{i}^{\frac{\beta}{\alpha}} )^T \bar{v}_{i} =   (\bar{q}_i +\bar{q}_{i}^{\frac{\beta}{\alpha}} )^T (\bar{v}_{i}-\hat{v}_{i}^{*}+\hat{v}_{i}^{*}) \\
= &  (\bar{q}_i +\bar{q}_{i}^{\frac{\beta}{\alpha}} )^T (\bar{v}_{i}-\hat{v}_{i}^{*}) -
 (\bar{q}_i +\bar{q}_{i}^{\frac{\beta}{\alpha}} )^T (\gamma_1\bar{q}_{i}^{\alpha}+\gamma_2\bar{q}_{i}^{\beta}).
   \end{align*}
   By Lemma \ref{ha1}, for $0<\alpha<1$, one has
  \begin{align}
&|\bar{v}_{i}-\hat{v}_{i}^{*}|=  |(\bar{v}_{i}^\frac{1}{\alpha})^\alpha-(\hat{v}_{i}^{*\frac{1}{\alpha}})^\alpha|  \nonumber\\
\leq & 2^{1-\alpha} |\bar{v}_{i}^\frac{1}{\alpha}-\hat{v}_{i}^{*\frac{1}{\alpha}}|^\alpha=2^{1-\alpha}|\varepsilon_{i}|^\alpha.  \label{vistar}
\end{align}
And, by Lemma \ref{qi},
  \begin{align}
 &(\bar{q}_i +\bar{q}_{i}^{\frac{\beta}{\alpha}} )^T  (\bar{v}_{i}-\hat{v}_{i}^{*})  \leq  2^{1-\alpha}
 |\bar{q}_i |^T    |\varepsilon_{i}|^\alpha
 + 2^{1-\alpha} |\bar{q}_{i}^{\frac{\beta}{\alpha}} |^T    |\varepsilon_{i}|^\alpha \nonumber\\
  \leq & 2^{1-\alpha}(\frac{1}{1+\alpha}\|\bar{q}_{i}^{\alpha+1}\|_1+\frac{\alpha}{1+\alpha}\|\varepsilon_{i}^{\alpha+1}\|_1 )\nonumber\\
  &+ 2^{1-\alpha}(\frac{\frac{\beta}{\alpha}}{\frac{\beta}{\alpha}+\alpha}\|\bar{q}_{i}^{\frac{\beta}{\alpha}+\alpha}\|_1+\frac{\alpha}{\frac{\beta}{\alpha}+\alpha}\|\varepsilon_{i}^{\frac{\beta}{\alpha}+\alpha}\|_1 ). \label{qTvi}
\end{align}
Using \eqref{vistar} and \eqref{qTvi}, one has
   \begin{align}
    \dot W_{1i} (\bar q_i)
\leq & -(\gamma_1-\frac{2^{1-\alpha}}{1+\alpha} ) \| \bar{q}_{i}^{\alpha+1}\|_1 -\gamma_2 \|\bar{q}_{i}^{\beta+1}\|_1 \nonumber\\
 & -(\gamma_1 -\frac{2^{1-\alpha}\frac{\beta}{\alpha}}{\frac{\beta}{\alpha}+\alpha})\| \bar{q}_{i}^{\frac{\beta}{\alpha}+\alpha}\|_1
 -\gamma_2 \| \bar{q}_{i}^{\frac{\beta}{\alpha}+\beta}\|_1  \nonumber \\
 & +\frac{2^{1-\alpha}\alpha}{\alpha+1} \| \varepsilon_{i}^{\alpha+1}\|_1
+\frac{2^{1-\alpha}\alpha}{\frac{\beta}{\alpha}+\alpha} \| \varepsilon_{i}^{\frac{\beta}{\alpha}+\alpha} \|_1. \label{dW1}
  \end{align}

Next, we define a vector function
  \begin{align*}
  d(\bar q_i, \bar v_i) =\int_{\hat{v}_{i}^{*}}^{\bar{v}_{i}}(s^{\frac{1}{\alpha}}-\hat{v}_{i}^{*\frac{1}{\alpha}})^{2-\alpha}ds.
\end{align*}
and hence
  \begin{align*}
W_{2i}(\bar q_i, \bar v_i) =   \|d(\bar q_i, \bar v_i)\|_1.
\end{align*}
Before the analysis on its derivative, we give the following calculation in order:
  \begin{align}
 &  \int_{\hat{v}_{i}^{*}}^{\bar{v}_{i}}(s^{\frac{1}{\alpha}}-\hat{v}_{i}^{*\frac{1}{\alpha}})^{1-\alpha}ds
  \leq
  \mbox{diag}(|\bar{v}_{i}-\hat{v}_{i}^{*}|) |\varepsilon_{i}|^{1-\alpha} , \nonumber \\
  &   |\bar{v}_{i}|^T  \mbox{diag}(|\bar{v}_{i}-\hat{v}_{i}^{*}|) |\varepsilon_{i}|^{1-\alpha} \leq 2^{1-\alpha}  |\bar{v}_{i}|^T  |\varepsilon_{i}| .  \label{intsds}
      \end{align}
Then, the derivative of $W_{2i}(\bar q_i, \bar v_i)$ along the trajectory of (\ref{dou}) satisfies,
 using  \eqref{intsds},
  \begin{align}
&\dot W_{2i}(\bar q_i, \bar v_i) \nonumber\\
=& (2-\alpha)
 \dot{\bar{q}}_{i}^T   \frac{-\partial \bar{v}_{i}^{*\frac{1}{\alpha}}}{\partial \bar{q}_{i}}
 \int_{\hat{v}_{i}^{*}}^{\bar{v}_{i}}(s^{\frac{1}{\alpha}}-\hat{v}_{i}^{*\frac{1}{\alpha}})^{1-\alpha}ds
    +(\varepsilon_{i}^{2-\alpha} )^T \dot{\bar{v}}_{i} \nonumber \\
\leq  & A_1 +A_2 \label{dW2}
    \end{align}
 for
 \begin{align*}
  &A_1 = (2-\alpha)2^{1-\alpha} |\bar{v}_{i}|^T \left|\frac{\partial (\hat{v}_{i}^{*\frac{1}{\alpha}})}{\partial \bar{q}_{i}}\right| |\varepsilon_{i}|\\
  &A_2 =(\varepsilon_{i}^{2-\alpha} )^T(u_{2i}-E\ddot{\eta}_i+Z_{i}).
 \end{align*}

By Lemma \ref{ha}, we can obtain
  \begin{align*}
   & \left|\frac{\partial (\hat{v}_{i}^{*\frac{1}{\alpha}})}{\partial \bar{q}_{i}}\right| \\
   =& \left| \mbox{diag}((\gamma_1\bar{q}_{i}^{\alpha}+\gamma_2\bar{q}_{i}^{\beta})^{\frac{1}{\alpha}-1})
   \mbox{diag}  (\gamma_1  \bar{q}_{i}^{\alpha-1}+\frac{\gamma_2 \beta}{\alpha}\bar{q}_{i}^{\beta-1})\right| \\
\leq & \gamma_1^{\frac{1}{\alpha}} I_m +\gamma_1^{\frac{1}{\alpha}-1}\frac{\gamma_2\beta}{\alpha}
 \mbox{diag} (|\bar{q}_{i}^{\beta-\alpha}|) \\
&+\gamma_1\gamma_2^{\frac{1}{\alpha}-1}  \mbox{diag}(|\bar{q}_{i}^{\frac{\beta}{\alpha}-\beta+\alpha-1} |)
+\gamma_2^{\frac{1}{\alpha}}  \frac{\beta}{\alpha}\mbox{diag} (|\bar{q}_{i}^{\frac{\beta}{\alpha}-1}|)
  \end{align*}
that implies
$
A_1\leq  (2-\alpha)2^{1-\alpha}   \sum_{j=1}^4
 \lambda_j |\bar{v}_{i}|^T\mbox{diag} (|\bar{q}_{i}^{p_j}|)|\varepsilon_{i}| .$

 To simplify the presentation, we introduce the following operator
   \begin{align*}
   \langle x, y \rangle_q  :=x \| \bar q_i^y \|_1,\; \langle x, y \rangle_\varepsilon  :=x \| \varepsilon_i^y \|_1.
    \end{align*}
Then,  using  Lemma \ref{qi} and a similar argument as  \eqref{vistar} gives
   \begin{align*}
  &\lambda |\bar{v}_{i}|^T \mbox{diag}( |\bar{q}_{i}^{p}|) |\varepsilon_{i}|\nonumber\\
  \leq &
  \lambda |\bar{v}_{i}-\hat{v}^*_{i}|^T  \mbox{diag}( |\bar{q}_{i}^{p}|) |\varepsilon_{i}| +\lambda |\hat{v}^*_{i}|^T \mbox{diag}( |\bar{q}_{i}^{p}|) |\varepsilon_{i}|  \nonumber\\
  \leq & \lambda 2^{1-\alpha}|\varepsilon_{i}^{\alpha+ 1}|^T  |\bar{q}_{i}|^{p}
  +\lambda|\gamma_1\bar{q}_{i}^{\alpha} +\gamma_2\bar{q}_{i}^{\beta} |\mbox{diag}( |\bar{q}_{i}^{p}|) |\varepsilon_{i}| \nonumber\\
    \leq & \lambda 2^{1-\alpha}|\varepsilon_{i}^{\alpha+ 1}|^T  |\bar{q}_{i}|^{p}
  +\lambda \gamma_1  |\varepsilon_{i}|^T |\bar{q}_{i}^{p+\alpha}| +\lambda \gamma_2  |\varepsilon_{i}|^T |\bar{q}_{i}^{p+\beta} | \nonumber\\
%
  %
%
  \leq &
  \langle l_1(p), p+\alpha+1 \rangle_q  +
\langle l_2(p), p+\beta+1 \rangle_q  \\
&+ \langle l_3(p,\lambda), p+\alpha+1 \rangle_\varepsilon +
\langle l_4(p,\lambda), p+\beta+1 \rangle_\varepsilon
 \end{align*}
and hence
   \begin{align}
A_1\leq &     \sum_{j=1}^4
 \langle \ell_1(p_j), p_j+\alpha+1 \rangle_q  +
\langle \ell_2(p_j), p_j+\beta+1 \rangle_q + \nonumber\\ &
\langle \ell_3(p_j,\lambda_j), p_j+\alpha+1 \rangle_\varepsilon +
\langle \ell_4(p_j,\lambda_j), p_j+\beta+1 \rangle_\varepsilon .  \label{A1}   \end{align}

By Lemma~\ref{lemma_fixobs}, $y_i(t)=0$ and hence  $SS\eta_i-\ddot{\eta}_i=0$ for $t\geq T_1^{*}$.
By Lemma~\ref{lemma_u1i}, one has $\zeta_i^TZ_i \leq 0$. Since $|\varepsilon_{ji}^{\alpha+1}|=\varepsilon_{ji}^{\alpha+1}$ and $|\varepsilon_{ji}^{\frac{\beta}{\alpha}+\beta}|=\varepsilon_{ji}^{\frac{\beta}{\alpha}+\beta}$, $j=1,\cdots,m$, for $\alpha$ and $\beta$ being two rational numbers of ratio of two odd integers and $\varepsilon_{ji}$ being the $j$-th entry of $\varepsilon_i$,
  \begin{align*}
(\varepsilon_{i}^{2-\alpha})^T   \varepsilon_{i}^{2\alpha-1} &= \sum_{j=1}^m\varepsilon_{ji}^{\alpha+1}=\sum_{j=1}^m|\varepsilon_{ji}^{\alpha+1}|=\| \varepsilon_{i}^{\alpha+1}\|_1 ,\\
 (\varepsilon_{i}^{2-\alpha})^T   \varepsilon_{i}^{\frac{\beta}{\alpha}+\beta+\alpha-2} &=\| \varepsilon_{i}^{\frac{\beta}{\alpha}+\beta}\|_1.
 \end{align*} Thus,
  \begin{align}
A_2 =&(\varepsilon_{i}^{2-\alpha} )^T( -k_1 \varepsilon_{i}^{2\alpha-1}-k_2 \varepsilon_{i}^{\frac{\beta}{\alpha}+\beta+\alpha-2} ) \nonumber\\
&+(\varepsilon_{i}^{2-\alpha} )^T(ESS\eta_i-E\ddot{\eta}_i)+ \zeta_i^TZ_i
 \nonumber\\
\leq & -k_1 \| \varepsilon_{i}^{\alpha+1}\|_1 -k_2   \|\varepsilon_{i}^{\frac{\beta}{\alpha}+\beta}\|_1. \label{A2}
\end{align}

Let
  \begin{align*}
  W_{i}(\bar q_i, \bar v_i) =
  W_{1i} (\bar q_i)  + W_{2i}(\bar q_i, \bar v_i).
\end{align*}
Under the conditions for $\gamma_1$, $\gamma_2$, $k_1$, and $k_2$,
there exist $\hat{\gamma}>0$ and $\hat{k}>0$ satisfying
  \begin{align*}
\gamma_1-\frac{2^{1-\alpha}}{1+\alpha} - \ell_1(p_1) \geq \hat \gamma +2L_1, \;
   \gamma_2 -\ell_2(p_4) \geq \hat \gamma +2L_1
\end{align*}
and
  \begin{align*}
k_1- \frac{2^{1-\alpha}\alpha}{\alpha+1} -\ell_3(p_1, \lambda_1)  \geq \hat k +4L_2,\;
k_2-\ell_4(p_4, \lambda_4) \geq \hat k +4L_2 .
\end{align*}
Then,  combining \eqref{dW1}, \eqref{dW2}, \eqref{A1}, and \eqref{A2} gives, for $t\geq T_1^{*}$,
%
%
\begin{align*}
  \dot W_{i}(\bar q_i, \bar v_i)  \leq  -B_1 - B_2,
\end{align*}
where
  \begin{align*}
B_1 =&
    \langle \hat \gamma +2L_1,  \alpha+1 \rangle_q    + \langle \hat \gamma +2L_1,  \frac{\beta}{\alpha}+\beta   \rangle_q \\
  &  -\langle L_1 , 2\beta-\alpha+1 \rangle_q   -\langle L_1,  \frac{\beta}{\alpha}-\beta+2\alpha  \rangle_q\\
  B_2 =&  \langle \hat k +4L_2  ,  \alpha+1 \rangle_\varepsilon  +\langle \hat k +4L_2 , \frac{\beta}{\alpha}+\beta \rangle_\varepsilon \\
  & - \langle L_2,  \frac{\beta}{\alpha}+\alpha \rangle_\varepsilon  - \langle  L_2  ,\beta+1 \rangle_\varepsilon \\
&  - \langle  L_2 ,2 \beta-\alpha+1 \rangle_\varepsilon   - \langle   L_2,\frac{\beta}{\alpha}-\beta+2\alpha \rangle_\varepsilon.
 \end{align*}
It is easy to verify the following inequalities
\begin{align*}
  &\frac{\beta}{\alpha}+\beta > 2\beta-\alpha+1>\alpha+1\\
  &\frac{\beta}{\alpha}+\beta > \frac{\beta}{\alpha}-\beta+2\alpha>\alpha+1\\
  &\frac{\beta}{\alpha}+\beta > \frac{\beta}{\alpha}+\alpha >\alpha+1\\
  &\frac{\beta}{\alpha}+\beta >\beta+1 >\alpha+1.
\end{align*}
Therefore,
  \begin{align*}
  &  \langle  2L_1,  \alpha+1 \rangle_q    + \langle 2L_1,  \frac{\beta}{\alpha}+\beta   \rangle_q \\
  \geq &\langle L_1 , 2\beta-\alpha+1 \rangle_q   + \langle L_1,  \frac{\beta}{\alpha}-\beta+2\alpha  \rangle_q
   \end{align*}
 that implies
$
B_1 \geq
    \langle \hat \gamma  ,  \alpha+1 \rangle_q    + \langle \hat \gamma ,  \frac{\beta}{\alpha}+\beta   \rangle_q.$
Similarly, one has
$
 B_2 \geq  \langle \hat k   ,  \alpha+1 \rangle_\varepsilon  +\langle \hat k   , \frac{\beta}{\alpha}+\beta \rangle_\varepsilon.$
The above two inequalities conclude
   \begin{align}
  \dot W_{i}(\bar q_i, \bar v_i)  \leq &  - \langle \hat \gamma  ,  \alpha+1 \rangle_q    - \langle \hat \gamma ,  \frac{\beta}{\alpha}+\beta   \rangle_q  \nonumber \\
  & - \langle \hat k   ,  \alpha+1 \rangle_\varepsilon  -\langle \hat k   , \frac{\beta}{\alpha}+\beta \rangle_\varepsilon. \label{dW}
\end{align}

Next, since
  \begin{align*}
\|d(\bar q_i, \bar v_i)\|_1 \leq \mbox{diag}( |\bar{v}_{i}-\hat{v}_{i}^{*}| ) |\varepsilon_i |^{2-\alpha}
\leq  2^{1-\alpha} \|\varepsilon_{i} \|^2,
  \end{align*}
 one has
  \begin{align*}
W_{i} (\bar q_i, \bar v_i) \leq \frac{1}{2}\|\bar{q}_{i}\|^2+\frac{1}{\beta/\alpha+1} \| \bar{q}_{i}^{\frac{\beta}{\alpha}+1}\|_1 +
2^{1-\alpha} \|\varepsilon_{i} \|^2.
 \end{align*}
Direct calculation, using Lemma~\ref{ha1}, gives
  \begin{align*}
  W_{i}^{\frac{\alpha+1}{2}} \leq &
 \langle \nu_1, \alpha+1 \rangle_q + \langle \nu_1, (\frac{\beta}{\alpha}+1)\frac{\alpha+1}{2}   \rangle_q +
  \langle \nu_1,  \alpha+1\rangle_\varepsilon
   \end{align*}
  and
    \begin{align*}
  W_{i}^{\frac{\frac{\beta}{\alpha}+\beta}{\frac{\beta}{\alpha}+1}} \leq &
 \langle \nu_2, \frac{2(\frac{\beta}{\alpha}+\beta )}{\frac{\beta}{\alpha}+1} \rangle_q + \langle \nu_2, \frac{\beta}{\alpha}+\beta    \rangle_q    +
\langle \nu_2, \frac{2(\frac{\beta}{\alpha}+\beta )}{\frac{\beta}{\alpha}+1} \rangle_\varepsilon
   \end{align*}
for some constants $\nu_1,\nu_2>0$. Again, It is easy to verify the following inequalities
       \begin{align*}
& \frac{\beta}{\alpha}+\beta > (\frac{\beta}{\alpha}+1)\frac{\alpha+1}{2}  >\alpha+1 \\
 & \frac{\beta}{\alpha}+\beta > \frac{2(\frac{\beta}{\alpha}+\beta)}{\frac{\beta}{\alpha}+1}>\alpha+1.
   \end{align*}
 Therefore,
 \begin{align}
  W_{i}^{\frac{\alpha+1}{2}} \leq &
 \langle 2 \nu_1, \alpha+1 \rangle_q + \langle \nu_1, \frac{\beta}{\alpha}+\beta   \rangle_q +
  \langle \nu_1,  \alpha+1\rangle_\varepsilon  \nonumber\\
  W_{i}^{\frac{\frac{\beta}{\alpha}+\beta}{\frac{\beta}{\alpha}+1}} \leq &
 \langle \nu_2,  \alpha+1 \rangle_q + \langle 2 \nu_2, \frac{\beta}{\alpha}+\beta    \rangle_q  \nonumber \\
 &  +
\langle \nu_2, \alpha+1 \rangle_\varepsilon
+ \langle \nu_2,   \frac{\beta}{\alpha}+\beta  \rangle_\varepsilon . \label{Wi}
   \end{align}

Comparing \eqref{dW} with \eqref{Wi}, one can conclude
     \begin{align*}
   & \dot{W}_{i}\leq  - \rho_1 W_{i}^{\frac{\alpha+1}{2}}
   -\rho_2 W_{i}^{\frac{\frac{\beta}{\alpha}+\beta}{\frac{\beta}{\alpha}+1}}
\end{align*}
for
$\rho_1 = \min \left\{ \frac{\hat \gamma} {2 \nu_1}, \frac{\hat \kappa} { \nu_1 } \right\} /2,\;
 \rho_2 = \min \left\{ \frac{\hat \gamma} {2\nu_2}, \frac{\hat \kappa} {\nu_2} \right\}/2.$
By Lemma~\ref{poly}, the equilibrium of  (\ref{dou}) is fixed-time stable.
In particular, there exists
  \begin{align*}
 T_2^*\leq  \frac{2}{\rho_1(1-\alpha)}+\frac{ \beta +\alpha}{\rho_2 \alpha ( \beta-1)}
  \end{align*}
  such that  \eqref{T2star} holds.
\end{Proof}

Finally,  based on Lemma~\ref{lemma_fixobs} and Lemma \ref{unii}, we can obtain the following theorem for the solvability
of the fixed-time synchronization problem with $T^* =T_1^*+T_2^*$.

\begin{Theorem}\label{ddf}
The fixed-time synchronization problem for
the multi-agent system composed of (\ref{1}) and (\ref{2}) under Assumption \ref{a1} is solvable by the observer
\eqref{ob} and the controller $\tau_i=\hat{M} (u_{1i}+u_{2i})$ of the form
\eqref{u1i} and \eqref{unfix} with all the parameters given in  Lemma~\ref{lemma_fixobs} and Lemma \ref{unii}.
\end{Theorem}

\begin{Remark}
Suppose the sub-controller $u_{1i}$ follows \eqref{u1i} with \eqref{zetai} but
the sub-controller \eqref{unfix}  for $u_{2i}$ reduces to the following finite-time controller,
by setting $k_2=0$ and $\gamma_2=0$,
  \begin{align} \label{un}
  u_{2i} &=-k_1 \varepsilon_{i}^{2\alpha-1} +ESS\eta_i, ~~
\varepsilon_{i}=\bar{v}_{i}^{\frac{1}{\alpha}} + \gamma_1^{\frac{1}{\alpha}} \bar{q}_{i},
\end{align}
where $\eta_i$ is governed by the finite-time observer \eqref{ob1}.
If  $\frac{1}{2}<\alpha<1$ is a ratio of two odd integers and $\gamma_1, k_1$ satisfy
  \begin{align*}
   &\gamma_1 > \frac{ 2^{1-\alpha}}{1+\alpha}+\frac{\alpha(2-\alpha) 2^{1-\alpha}}{1+\alpha} \\
   &k_1 > \gamma_1^{1+1/\alpha} \left(\frac{2^{1-\alpha}\alpha}{1+\alpha}+
\frac{(2-\alpha) 2^{1-\alpha}}{\gamma_1}(2^{1-\alpha}+\frac{\gamma_1}{1+\alpha}) \right),\nonumber
  \end{align*}
then the equilibrium of the closed-loop system composed of (\ref{dou})
at the origin is finite-time stable. In particular, there exists
a finite settling-time function
$T_{2i}(\bar q_i(T_1(\eta(0))) ,\bar v_i(T_1(\eta(0))) ) \geq 0$ such that
\begin{align}
& \lim_{t\rightarrow T_1+ T_{2i}  } [\bar q_i(t), \bar v_i(t)] =0,\\ \label{T2}
&  [\bar q_i(t), \bar v_i(t)]=0,\; t\geq T_1+T_{2i} , \; \forall  \bar q_i(T_1), \bar v_i(T_1) \in \mathbb{R}^{m}.\nonumber
\end{align}
As a result,  the finite-time synchronization problem for
the multi-agent system composed of (\ref{1}) and (\ref{2})
under Assumption \ref{a1} is solvable by the observer
\eqref{ob1} and the controller $\tau_i=\hat{M} (u_{1i}+u_{2i})$ of the form
\eqref{u1i} and \eqref{un}. The proof can similarly follow that of  Lemma~\ref{unii} and is thus omitted.
 \end{Remark}

\section{An Example}

\begin{figure}[t]
  \begin{center}
  \includegraphics[width=0.6\hsize, bb=0 10 160 70]{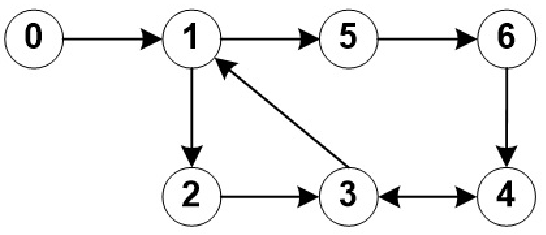}\\
  \end{center}
  \caption{Illustration of the communication network topology.}\label{g}
\end{figure}

Consider a group of six  robotic manipulators given by (\ref{1}) where $q_i=[q_{1i},q_{2i}]^T\in\mathbb{R}^2$ and
\begin{equation*}\label{pa}
  \begin{aligned}
 & M_i(q_i)=\left[
        \begin{array}{cc}
          \theta_{1i}+\theta_{2i}+2\theta_{3i}\cos(q_{2i}) & \theta_{2i}+\theta_{3i}\cos(q_{2i}) \\
          \theta_{2i}+\theta_{3i}\cos(q_{2i}) &\theta_{4i} \\
        \end{array}
      \right]\\
&C_{i}(q_i,\dot{q}_i)\dot{q}_i=\left[
        \begin{array}{cc}
          -\theta_{3i}\sin(q_{2i})\dot{q}_{1i}^2-2\theta_{3i}\sin(q_{2i})\dot{q}_{1i}\dot{q}_{2i} \\
          \theta_{3i}\sin(q_{2i})\dot{q}_{2i}^2 \\
        \end{array}
      \right]  \\
&G_i(q_i)=\left[
        \begin{array}{cc}
          \theta_{5i}g\cos(q_{1i})+\theta_{6i}g\cos(q_{1i}+q_{2i}) \\
          \theta_{6i}g\cos(q_{1i}+q_{2i}) \\
        \end{array}
      \right]\\
  \end{aligned}
\end{equation*}
for $i=1, \cdots, 6$. In the equations, $\theta_{ji}$, $j=1,\cdots,6$, $i=1,\cdots,6$, represent unknown parameters.
The leader system is given by (\ref{2}) with $S=\left[
                                          \begin{array}{ccc}
                                            0&1 \\
                                            -1&0\\
                                          \end{array}
                                        \right]
$ and $E=I_2$.
The information flow among all the subsystems and the leader is described by the digraph in Fig.~\ref{g}, which contains a spanning tree with node 0 as the root, satisfying Assumption \ref{a1}.
Let $D=8I_6$. Then $DH+H^TD\geq 2I_N$.
\begin{figure}[t]
  \begin{center} \includegraphics[width=\hsize]{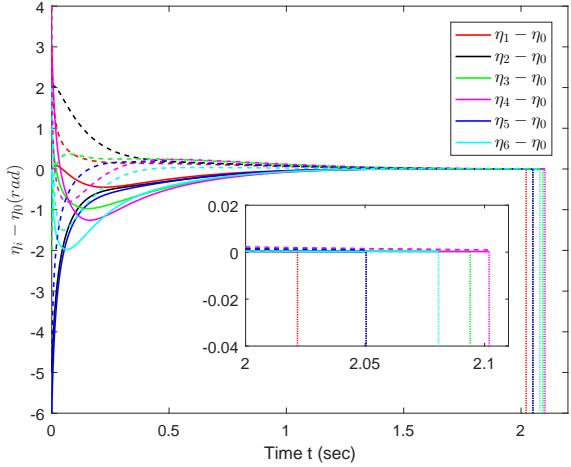}\\
  \end{center}
  \caption{ Profile of the estimation errors $\eta_i-\eta_0$, $i=1,\cdots,6$,
under the  fixed-time observer.}\label{gg1}
\end{figure}

By  Lemma \ref{lemma_fixobs}, let
$c_1=8.4,~c_2=1,~c_3=1,~a=3/5,~b=3$. Now we can construct the fixed-time observer (\ref{ob}) whose performance
is shown in Fig. \ref{gg1}. It is observed that the estimation errors $\eta_i-\eta_0$, $i=1,\cdots,6$, approach zero  at the
time instants marked by the vertical lines.
In the simulation, the error tolerance of numerical calculation is set as $10^{-3}$ that is used as the criterion of approaching zero.

\begin{figure}[t]
\centering
  \includegraphics[width=\hsize]{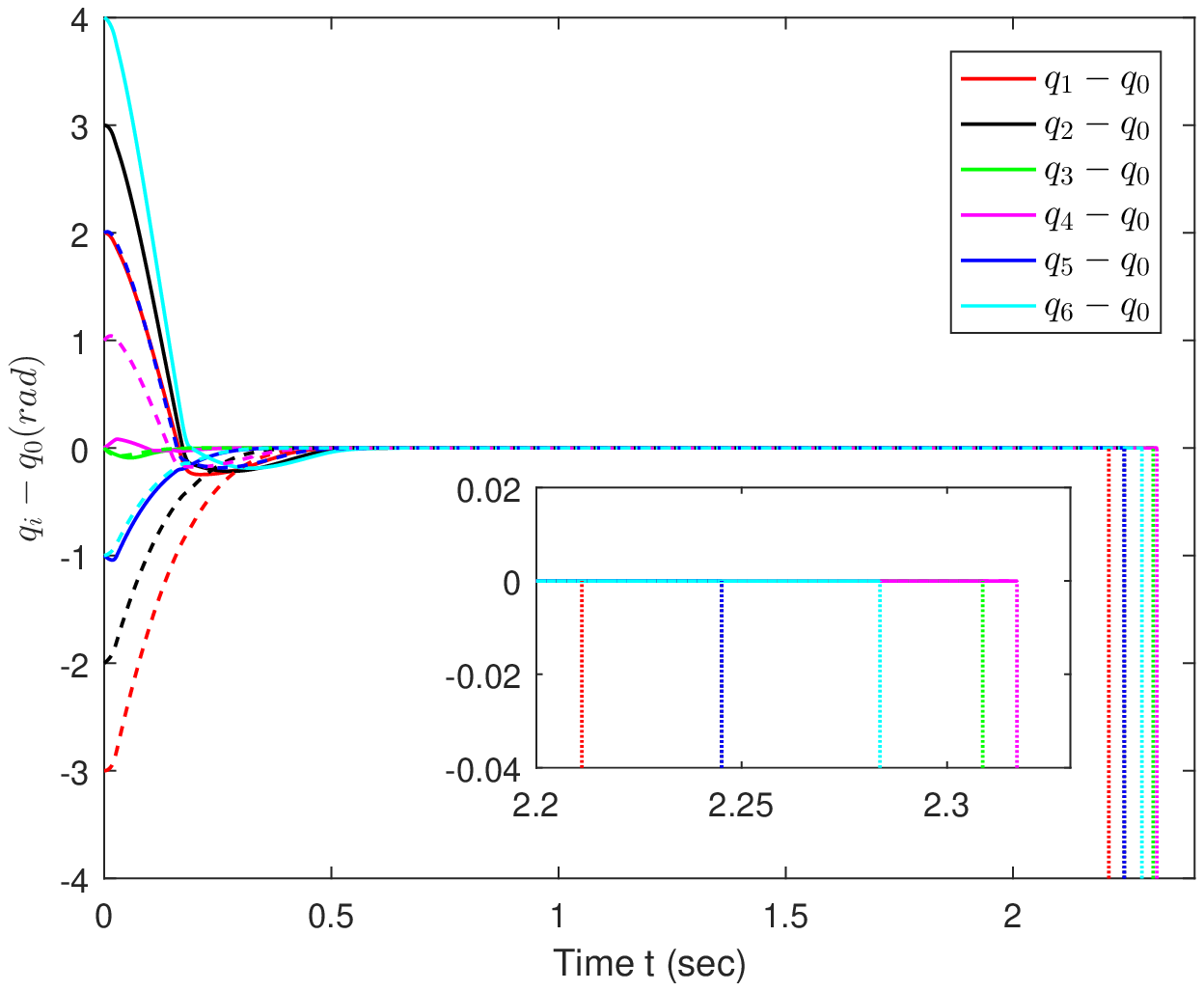}
    \includegraphics[width=\hsize]{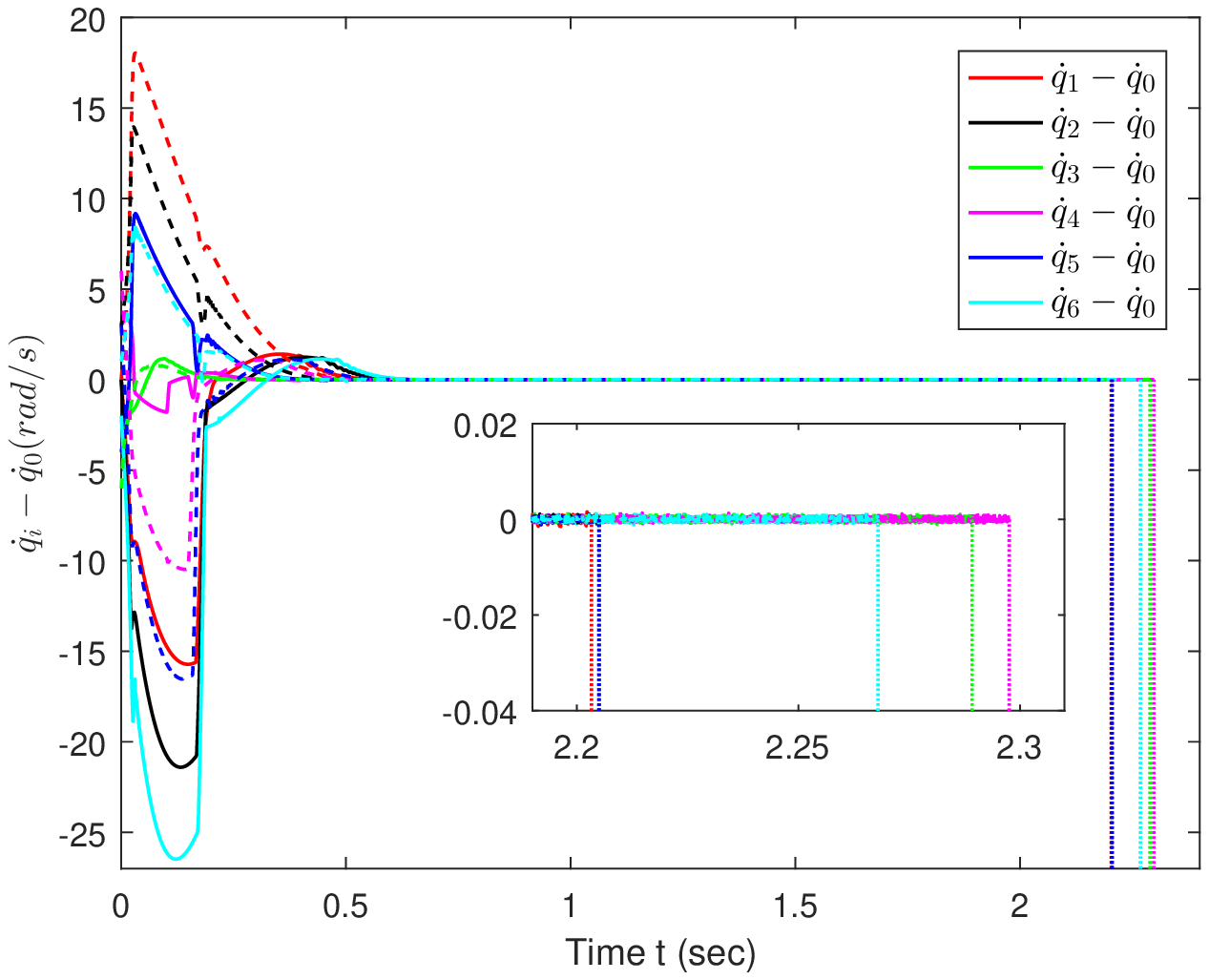}
  \caption{Profile of the synchronization errors $q_i-q_0$ and $\dot{q}_i-\dot{q}_0$, $i=1,\cdots,6$, under the fixed-time controller.}\label{g3}
\end{figure}

Next, we apply the observer (\ref{ob}) to solve the fixed-time control problem of Euler-Lagrange systems and design the fixed-time control law $\tau_i=\hat{M}(u_{1i}+u_{2i})$ where $u_{1i}$ is given by (\ref{u1i}) and $u_{2i}$ is given by (\ref{unfix}). Although we we do not know the exact value of $M_i(q_i)$, $C_i(q_i,\dot{q}_i)$ and $G_i(q_i)$,
it is assumed that the unknown parameters in the following ranges
$ \theta_{1i} \in [6,8]$, $ \theta_{2i} \in [0.8,1]$, $ \theta_{3i} \in [1,1.4]$, $\theta_{5i}\in[1.5,2]$, and $\theta_{6i}\in[1,1.3]$.
Simple calculation verifies that
the properties in \eqref{km} are satisfied for
$k_{\underline{m}}^{-1}=0.3,~k_{\overline{m}}^{-1}=0.08,~k_c=3$, and $k_g=50$.
We select the  parameters in (\ref{u1i}) and  (\ref{unfix}) as
$\kappa=3,~\epsilon=11/19,~\gamma_1=10,~\gamma_2=10,~k_1=20,~k_2=15,~\alpha={7}/{9},~\beta={9}/{7}$.
For the purpose of simulation, we provide the values for uncertain parameters
$\theta_{1i}=7,~\theta_{2i}=0.96,~\theta_{3i}=1.2,~\theta_{4i}=5.96,~\theta_{5i}=2$, and $\theta_{6i}=1.2$.
 The simulation is conducted with arbitrarily selected initial conditions.
Fig.~\ref{g3}   shows $q_i$, $\dot{q}_i$ respectively converge to $q_0$, $\dot{q}_0$ in fixed time instants.


\section{Conclusion}
This paper has proposed the fixed-time robust control  design for the consensus problem of networked Euler-Lagrange systems based on a distributed observer, which is capable of estimating the desired trajectory of the leader in a fixed time under a directed graph. The heterogeneous uncertain Euler-Lagrange systems are converted into second-order systems by a partial design of the control law, and then the backstepping procedure for second-order systems are utilized to accomplished  the fixed-time control design.

\bibliographystyle{IEEEtran}

\bibliography{finite}

\begin{thebibliography}{10}
\providecommand{\url}[1]{#1}
\csname url@samestyle\endcsname
\providecommand{\newblock}{\relax}
\providecommand{\bibinfo}[2]{#2}
\providecommand{\BIBentrySTDinterwordspacing}{\spaceskip=0pt\relax}
\providecommand{\BIBentryALTinterwordstretchfactor}{4}
\providecommand{\BIBentryALTinterwordspacing}{\spaceskip=\fontdimen2\font plus
\BIBentryALTinterwordstretchfactor\fontdimen3\font minus
  \fontdimen4\font\relax}
\providecommand{\BIBforeignlanguage}[2]{{%
\expandafter\ifx\csname l@#1\endcsname\relax
\typeout{** WARNING: IEEEtran.bst: No hyphenation pattern has been}%
\typeout{** loaded for the language `#1'. Using the pattern for}%
\typeout{** the default language instead.}%
\else
\language=\csname l@#1\endcsname
\fi
#2}}
\providecommand{\BIBdecl}{\relax}
\BIBdecl

\bibitem{bhat2000}
S.~P. Bhat and D.~S. Bernstein, ``Continuous finite-time stabilization of the
  translational and rotational double integrators,'' \emph{IEEE Transactions on
  Automatic Control}, vol.~43, no.~5, pp. 678--682, 1998.

\bibitem{xiao2009}
F.~Xiao, L.~Wang, J.~Chen, and Y.~Gao, ``Finite-time formation control for
  multi-agent systems,'' \emph{Automatica}, vol.~45, no.~11, pp. 2605--2611,
  2009.

\bibitem{fu2016}
J.~Fu and J.~Wang, ``Observer-based finite-time coordinated tracking for
  general linear multi-agent systems,'' \emph{Automatica}, vol.~66, pp.
  231--237, 2016.

\bibitem{cao2014}
Y.~Cao and W.~Ren, ``Finite-time consensus for multi-agent networks with
  unknown inherent nonlinear dynamics,'' \emph{Automatica}, vol.~50, no.~10,
  pp. 2648--2656, 2014.

\bibitem{du2020}
H.~Du, G.~Wen, D.~Wu, Y.~Cheng, and J.~L{\"u}, ``Distributed fixed-time
  consensus for nonlinear heterogeneous multi-agent systems,''
  \emph{Automatica}, vol. 113, p. 108797, 2020.

\bibitem{zuo2015}
Z.~Zuo, ``Nonsingular fixed-time consensus tracking for second-order
  multi-agent networks,'' \emph{Automatica}, vol.~54, pp. 305--309, 2015.

\bibitem{zuo2017}
Z.~Zuo, B.~Tian, M.~Defoort, and Z.~Ding, ``Fixed-time consensus tracking for
  multiagent systems with high-order integrator dynamics,'' \emph{IEEE
  Transactions on Automatic Control}, vol.~63, no.~2, pp. 563--570, 2018.

\bibitem{hong2002}
Y.~Hong, Y.~Xu, and J.~Huang, ``Finite-time control for robot manipulators,''
  \emph{Systems \& Control Letters}, vol.~46, no.~4, pp. 243--253, 2002.

\bibitem{hong2006}
Y.~Hong, J.~Wang, and D.~Cheng, ``Adaptive finite-time control of nonlinear
  systems with parametric uncertainty,'' \emph{IEEE Transactions on Automatic
  Control}, vol.~51, no.~5, pp. 858--862, 2006.

\bibitem{lin1}
X.~Huang, W.~Lin, and B.~Yang, ``Global finite-time stabilization of a class of
  uncertain nonlinear systems,'' \emph{Automatica}, vol.~41, pp. 881--888,
  2005.

\bibitem{yu2005}
S.~Yu, X.~Yu, B.~Shirinzadeh, and Z.~Man, ``Continuous finite-time control for
  robotic manipulators with terminal sliding mode,'' \emph{Automatica},
  vol.~41, no.~11, pp. 1957--1964, 2005.

\bibitem{zhao2009}
D.~Zhao, S.~Li, and F.~Gao, ``A new terminal sliding mode control for robotic
  manipulators,'' \emph{International Journal of control}, vol.~82, no.~10, pp.
  1804--1813, 2009.

\bibitem{galicki2015}
M.~Galicki, ``Finite-time control of robotic manipulators,'' \emph{Automatica},
  vol.~51, pp. 49--54, 2015.

\bibitem{huang2015}
J.~Huang, C.~Wen, W.~Wang, and Y.-D. Song, ``Adaptive finite-time consensus
  control of a group of uncertain nonlinear mechanical systems,''
  \emph{Automatica}, vol.~51, pp. 292--301, 2015.

\bibitem{hu2019finite}
H.-X. Hu, G.~Wen, W.~Yu, J.~Cao, and T.~Huang, ``Finite-time coordination
  behavior of multiple {E}uler--{L}agrange systems in cooperation-competition
  networks,'' \emph{IEEE transactions on cybernetics}, vol.~49, no.~8, pp.
  2967--2979, 2019.

\bibitem{polyakov2011}
A.~Polyakov, ``Nonlinear feedback design for fixed-time stabilization of linear
  control systems,'' \emph{IEEE Transactions on Automatic Control}, vol.~57,
  no.~8, pp. 2106--2110, 2012.

\bibitem{su2015}
Y.~Su, ``Cooperative global output regulation of second-order nonlinear
  multi-agent systems with unknown control direction,'' \emph{IEEE Transactions
  on Automatic Control}, vol.~60, no.~12, pp. 3275--3280, 2015.

\bibitem{horn1994}
R.~A. Horn, R.~A. Horn, and C.~R. Johnson, \emph{Topics in Matrix
  Analysis}.\hskip 1em plus 0.5em minus 0.4em\relax Cambridge University Press,
  1994.

\bibitem{hardy1952}
G.~H. Hardy, J.~E. Littlewood, and G.~P{\'o}lya, \emph{Inequalities}.\hskip 1em
  plus 0.5em minus 0.4em\relax Cambridge University Press, 1952.

\bibitem{qian2001}
C.~Qian and W.~Lin, ``A continuous feedback approach to global strong
  stabilization of nonlinear systems,'' \emph{IEEE Transactions on Automatic
  Control}, vol.~46, no.~7, pp. 1061--1079, 2001.

\bibitem{song2017}
Y.~Song, Y.~Wang, J.~Holloway, and M.~Krstic, ``Time-varying feedback for
  regulation of normal-form nonlinear systems in prescribed finite time,''
  \emph{Automatica}, vol.~83, pp. 243--251, 2017.

\end{thebibliography}

\vfill

\end{document}